\documentclass{birkjour}
 \newtheorem{thm}{Theorem}[section]
 \newtheorem{cor}[thm]{Corollary}
 \newtheorem{lem}[thm]{Lemma}
 \newtheorem{prop}[thm]{Proposition}
 \theoremstyle{definition}
 
 \theoremstyle{remark}
 \newtheorem{rem}[thm]{Remark}
 
 \numberwithin{equation}{section}
\usepackage{epsfig}
\usepackage{amssymb, graphicx, amsmath, amsfonts}
\usepackage{multicol, hyperref, caption, subcaption, enumerate, longtable, enumitem}
\usepackage{array}
\usepackage{csquotes, rsfso, enumitem, multirow}

\newcolumntype{C}[1]{>{\centering\let\newline\\\arraybackslash\hspace{0pt}}m{#1}}

\begin{document}

\title{Forcing nonperiodic tilings with one tile \\ using a seed}



\author{Bernhard Klaassen} 

\begin{flushleft}
 Preprint, to appear 2021 in \emph{European Journal of Combinatorics} 
\end{flushleft}
 
\address{Fraunhofer Institute SCAI, Sankt Augustin\br
Germany}
\email{bernhard.klaassen@scai.fraunhofer.de}



\date{September 19, 2019}

\begin{abstract}
The so-called ``einstein problem'' (a pun playing with the famous scientist's name and the German term ``ein Stein'' for ``one stone'') asks for a simply connected prototile only allowing nonperiodic tilings without need of any matching rule. So far, researchers come only close to this demand by defining decorated prototiles forcing nonperiodicity of any generated tiling using matching rules. In this paper a class of spiral tilings (and one non-spiral example) is linked to a weaker form of the einstein problem where one or several seed tiles are used. Furthermore, the classical types of matching rules are listed and some new types are discussed. 
\end{abstract}

\maketitle

\section{Introduction}
\label{intro}
In 2011 studies on the ``einstein problem'' (in German ``ein Stein'' $=$ ``one stone'') came to a remarkable result, when
Socolar and Taylor~\cite{Aaht} presented their decorated single prototile forcing nonperiodicity of any generated
tiling by two matching rules. Other attempts were made also with one tile and two matching rules~\cite{dendrites}.
The ideal  ``einstein'' would be a single tile (a topological disk) only allowing nonperiodic tilings without need of
any matching rule. As long as such a perfect solution is unknown, it is an interesting task to come as close to it
as possible. So, in this paper, we are going to reduce the number of rules for certain prototiles to a single one.
(The difference to other approaches is that in our investigation the decoration on the tiles will generate a
continuous, self-avoiding curve.) Also the rules themselves will be grouped and discussed.

First, let us define the necessary terms: A plane \emph{tiling} $\mathcal{T}$ (or \emph{tessellation}) is a countable family of sets, called \emph{tiles}, that cover the plane without gaps or overlaps of non-zero area. If the intersection of three or more tiles is nonempty, then this intersection is called a \emph{vertex} of $\mathcal{T}$. An \emph{edge} is (part of) the intersection of two tiles that connects two vertices. In this research, all tiles considered are closed topological disks and uniformly bounded. Here we only consider $k$-hedral tilings, which means that each tiling has only finitely many congruence classes. The $k$ tiles, each representing a different class, are the \emph{prototiles}. With $k=1$, a tiling is called \emph{monohedral}.

To study a plane tiling $\mathcal{T}$, one usually looks at the different isometries that leave $\mathcal{T}$ unchanged. These isometries are called \emph{symmetries of $\mathcal{T}$} and the group of all symmetries of $\mathcal{T}$ is called the \emph{symmetry group of $\mathcal{T}$} denoted by $\mathcal{S(T)}$.

\section{A prototile and some of its properties}

\begin{figure}[!htb]
\centering
\includegraphics[width=0.55\columnwidth]{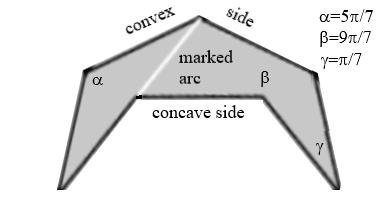}%
\caption{Tile $T$ with a marked arc (white line)}%
\label{Tile_T}%
\end{figure}
Consider tile $T$ as shown in Fig. \ref{Tile_T} with a marked arc represented by the white line. $T$ is constructed from of a regular heptagon and -- without the decoration -- was first discussed by Gr\"unbaum and Shephard
in their renowned standard work~\cite{TP}. In another wordplay they named such nonconvex tiles as
``versatiles''~\cite{versatiles} because of their multiple use in generating beautiful tilings.
(See~\cite{tilingsearch} for a large and nicely colored web collection not only with this type of tiles.) During
this paper, we distinguish two ``sides'' of the tile's boundary separated by the vertices with acute angles: the side
with inner angles greater than $\pi$ will be called ``concave side'' (the other side consequently is called
``convex side'').

We can now state a single rule for each newly placed tile during the process of building a tiling using copies of this prototile.  (As usual in the world of tilings, flipped copies of the prototile should be allowed, since mirror symmetry is an isometry.) Here we assume that a single tile was placed anywhere in the plane as starting set of the tiling or ``seed'' and the other tiles are added one by one in a sequence, obeying the

\noindent
\begin{itemize}[label={}]
\item \textbf{Matching rule R1}: Each additionally placed tile continuously connects with its marked arc to a marked arc of an already placed tile.
\end{itemize}

Such a combination of a seed and a classical matching rule is not a new concept in the context of nonperiodic tilings (e.g. in \cite{dendrites} or \cite{grow}). 

\begin{figure}[!htb]
\centering
\includegraphics[width=0.7\columnwidth]{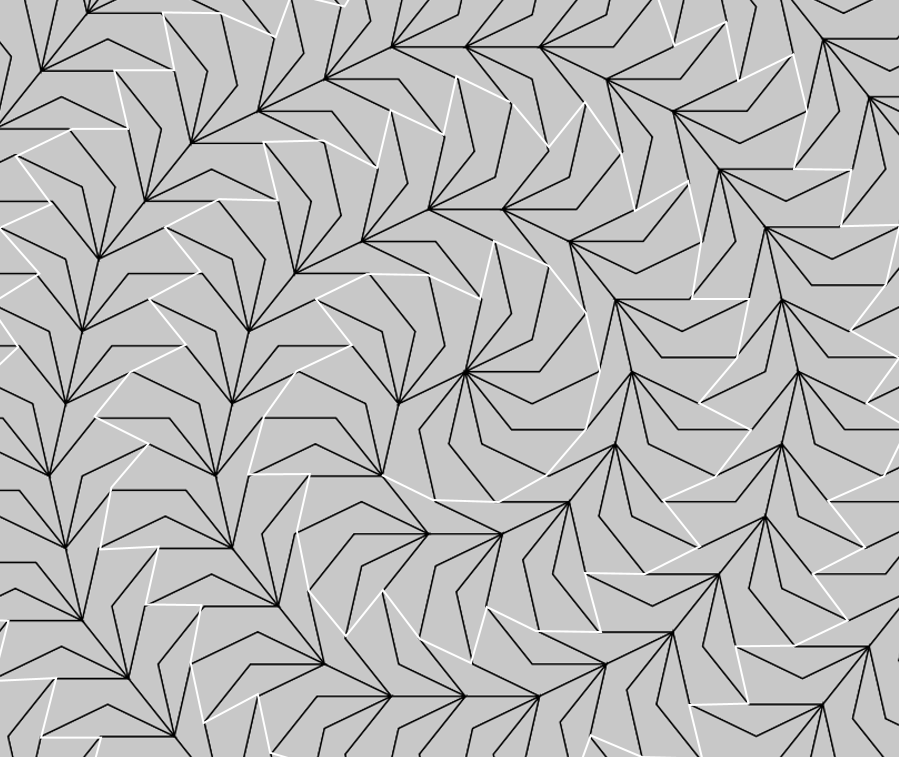}%
\caption{A tiling with prototile  $T$  obeying rule R1}%
\label{spiral}%
\end{figure}

By demonstrating a tiling of the plane (Fig. \ref{spiral}) with $T$ as prototile and obviously fulfilling rule R1, we can make sure that at least one tiling of this kind is possible. The union of marked arcs is represented by the white line. We'll see that tilings with prototile  $T$  obeying rule R1 must be nonperiodic. (Throughout this paper for simplicity sometimes the term ``arc(s)" will be used in the same meaning as ``marked arc(s)".)

Although the tiling sketched in Fig. \ref{spiral} is known from the literature, in the context of nonperiodic tilings one should give a proof that this really generates a tiling of the full plane.

\begin{figure}[!htb]
\centering
\includegraphics[width=0.7\columnwidth]{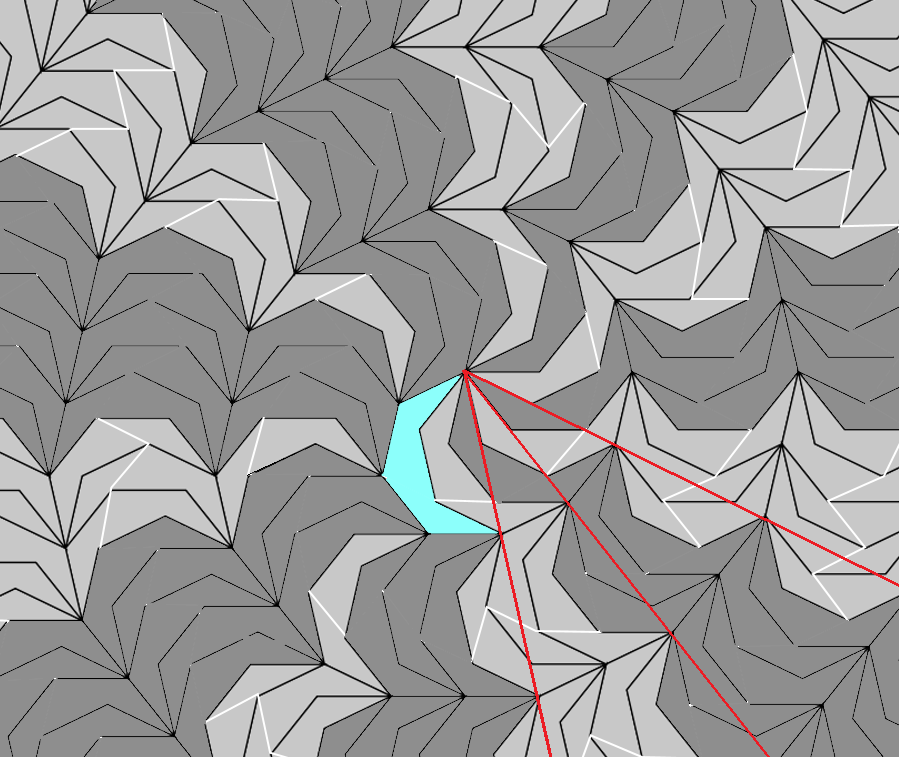}%
\caption{Partition of the spiral tiling into sectors}%
\label{sectors}%
\end{figure}

 This can be shown by partitioning the tiling into 14 sectors, see Fig. \ref{sectors}.
Note that each of these sectors has principally the same structure: It is a sequence of rows each of which consists of two tiles more than its predecessor. So, for most of the sectors the number of tiles per row is denoted by (1, 3, 5,\ldots) and three of them -- at the left half -- are characterized by (3, 5, 7,\ldots), i.e.,~the leading single tile was skipped, but apart from this first tile all sectors are congruent.  We can observe that each boundary between any two sectors is an unlimited curve periodically meeting a straight line (see the lines starting from the center in Fig. \ref{sectors}). Therefore, all sectors fit together without leaving any gaps and they grow beyond any limit covering the whole plane. This makes the union of the 14 sectors a tiling of the full plane.

For those readers who are not so familiar with spiral tilings it should be noted that all known spiral tilings
(see~\cite{tilingsearch}) can be partitioned into sectors in a more or less similar way.

Now, a simple observation can be made in form of

\begin{lem} \label{lem:curve}
The union of all marked arcs in a tiling generated by rule R1 with prototile $T$ is a continuous and self-avoiding curve. (See Figure $\ref{spiral}$.)
\end{lem}
\begin{proof}  By a simple induction argument starting with a single tile, we see that by adding arcs in a connected way to the curve of arcs obeying R1, this curve remains connected. To prove the self-avoidance, consider the inner angles at the endpoints of the marked arc ($5\pi/7$ and $9\pi/7$). Both inner angles are larger than $2\pi/3$. Therefore, it is impossible that three or more arcs meet at the same point, which would be necessary for a self-intersection of the resulting curve.
\end{proof}

\begin{rem} \label{rem:1}
 Observe that the statement of Lemma \ref{lem:curve} - ``the union of marked arcs is a continuous and self-avoiding curve'' - is not only concluded from R1, it is even equivalent to R1 as long as we assume a single tile as starting set for the tiling. Note that throughout this paper we assume such a decoration on each single tile to be of finite length.
\end{rem}

\begin{rem} \label{rem:2}
In the context of matching rules, it is often pointed out in the literature that the rules should be locally
checkable (preferably between neighboring tiles). In our case with rule R1 we also can perform a local check between
adjacent tiles beginning with the tile that we called ``starting set'' joined with the tile placed second. However, a certain order has to be kept in form of a sequence of local checks. A similar rule where small subtrees (``dendrites'') are added to a connected tree is used in~\cite{dendrites}.
\end{rem} 

Usually in the context of plane tilings, ``periodic'' means that two linear independent translation vectors exist, spanning a fundamental region. One could conclude that nonperiodic means the opposite of periodic and some authors use the term in this way. However, in the field of aperiodic tile sets, we should exclude also cases with only one translational symmetry. So, in this study ``nonperiodic'' means the lack of any translational symmetry. With our simple matching rule and the decorated prototile $T$, we can force nonperiodicity which is formally shown in 

\begin{prop} \label{prop:1}
Given a monohedral tiling with marked tiles generating a continuous self-avoiding curve of these arcs reaching each tile, such a tiling cannot have a translational symmetry.
\end{prop}
\begin{proof}
For the special case of a one-armed spiral shown in Fig. \ref{spiral}, the proof is self-evident, since the point where the curve (the union of arcs) begins, is unique for the whole tiling and cannot be found anywhere else. However, for other cases, where this curve can be extended in both directions, we need a more general proof.
Assume such a tiling $\mathcal{T}$  with a curve (through every tile) generated from the tiles' decoration and $\mathcal{S(T)}$ containing a translation by a nonzero vector $\tau$. For simplicity, we assume $\tau$ to be horizontal (see Fig. \ref{vector} for a visualization). 

\begin{figure}[!htb]
\centering
\includegraphics[width=0.75\columnwidth]{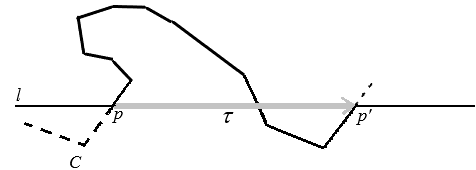}%
\caption{Assumed translation symmetry by vector $\tau$}%
\label{vector}%
\end{figure}

Let $p$ be a point where one of the marked arcs begins. Translation by vector $\tau$ yields point $p'$. We know that curve $C$ (composed of the arcs) must meet every marked arc in this tiling and can be considered as a sequence, i.e.,~starting from zero at the seed tile, each tile can be assigned a number from $\mathbb{Z}$ in the canonical order.  Hence, there are only finitely many arcs on the path between $p$ and $p'$, each of which with finite length. Therefore, the part of $C$ between $p$ and $p'$ (called $C_{p,p'}$ for short) as a continuous curve must be limited, and the maximal distance of $C_{p,p'}$ from the line  through $p$ and $p'$ (let it be named $l$) is also finite. Since the tiling has translational symmetry, curve $C$ as a periodic, continuous and self-avoiding  curve can be fully generated as union of periodical copies of $C_{p,p'}$ by infinite translation in direction  $\pm\tau$.
Therefore, the said maximal distance from line $l$ holds for the whole curve.   However, since the tiles are bounded and all must be met by $C$, there must be tiles in arbitrarily large distance from line $l$ traversed by $C$, which yields a contradiction.
 \end{proof}

\begin{cor} \label{cor}
 A tiling built from copies of tile $T$ and obeying rule R1 fulfills the assumptions of Proposition \ref{prop:1} and hence cannot be periodic.
 \end{cor}
\begin{proof}
The self-avoidance for the curve of marked arcs follows from Lemma \ref{lem:curve}.
\end{proof}

By demonstrating a tiling of the plane (Fig. \ref{spiral}) with $T$ as prototile and obviously fulfilling rule R1, we can make sure that at least one tiling of this kind is possible. And with Proposition \ref{prop:1} we know that this example and any other tiling generated by $T$ and rule R1 cannot be periodic.

\begin{figure}[!htb]
\centering
\includegraphics[width=0.55\columnwidth]{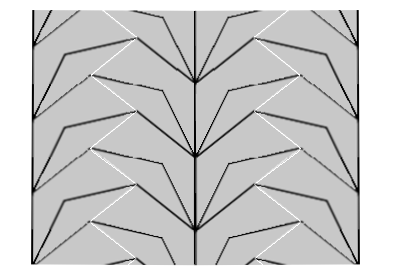}%
\caption{Periodic tiling with tile $T$ violating rule R1}%
\label{repeat}%
\end{figure}

On the other hand, it is not difficult to form a repeating structure with copies of tile $T$, e.g.,~see Fig. \ref{repeat}. Does this contradict our corollary? We will see that it does not: Also the curves formed by the unions of the marked arcs must be repeating, i.e.,~they are separated into disconnected parts according to the ``columns'' shown in the figure.

 Why is rule R1 violated here? While generating the repeating structure starting with a single seed tile, we can only build one column. The first tile belonging to a new column must violate rule R1.

It is interesting to note that Gr\"unbaum and Shephard used a similar (only slightly more general) rule in an early
attempt to define the term ``spiral tiling'' for monohedral tilings (see the exercise section of chapter 9.5.
in~\cite{TP}). They proposed that in the prototile a certain arc should be marked and if the union of all marked
arcs (from all tiles) can be partitioned into finitely many unlimited curves, then the tiling should be called
``spiral tiling''. The authors of ``Tilings and Patterns'' themselves found that this short and elegant definition had
the drawback that it was not suitable for all tilings with spiral character (and that it also worked for unwanted
cases). Therefore, in a 2017 paper~\cite{HTDST} (with a refinement in~\cite{detect}) a new definition was suggested,
which had to be a bit longer and more technical, but was applicable to all known examples of tilings with spiral
construction.

However, here we might apply the mentioned definition of Gr\"unbaum and Shephard to this tiling and find that rule R1 forces the tiling to be spiral and thus cannot be periodic.

\section{More examples}

There are several other interesting monohedral tilings which are also suitable for this approach. 
\begin{figure}[!htb]
\centering
\includegraphics[width=0.7\columnwidth]{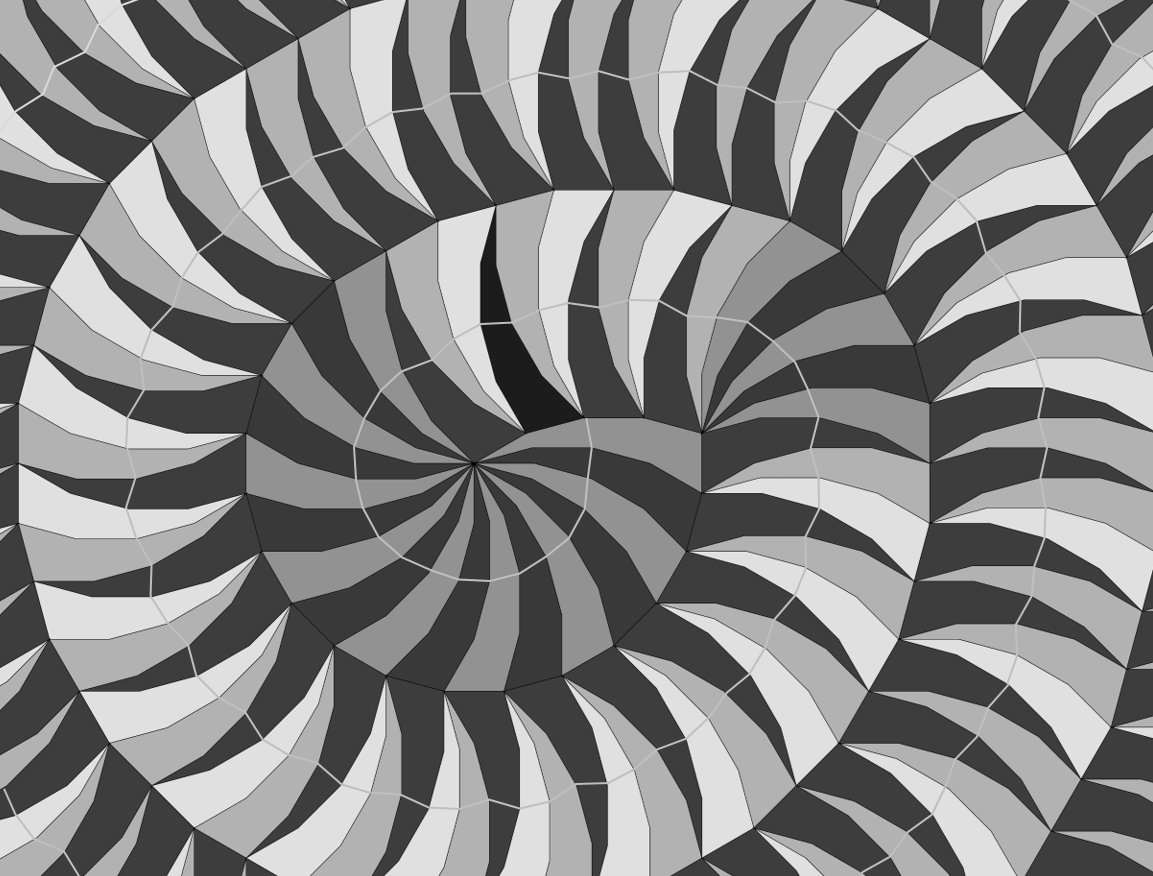}%
\caption{A tiling generated by another ``versatile''}%
\label{snake}%
\end{figure}
Our next example is again a ``versatile" from an old Grünbaum paper \cite{versatiles} shown in Fig. \ref{snake}. Again marked arcs are depicted crossing the tile (represented by a grey line).  The resulting nonperiodic tiling has the form of a nicely curled snake.

From Corollary \ref{cor} we know that for each new example we must check if the inner angles at the points where the arc meets the tile's boundary are larger than $2\pi/3$. This is the case, so Lemma \ref{lem:curve} and the proposition remains valid. 

\begin{figure}[!htb]
\centering
\includegraphics[width=0.65\columnwidth]{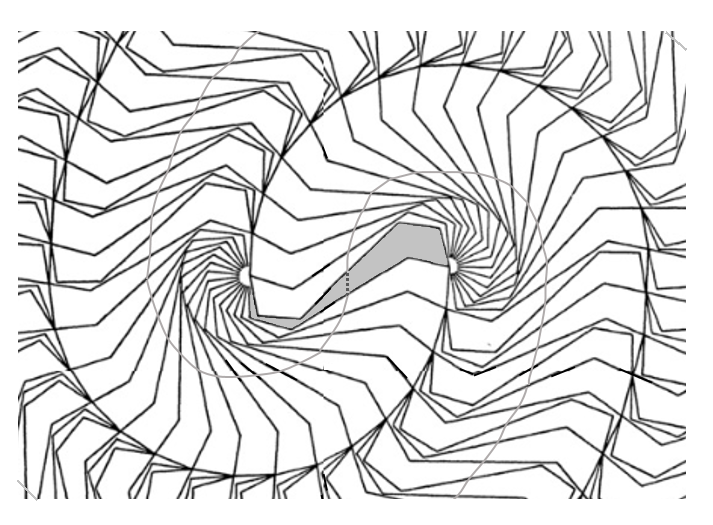}%
\caption{The Voderberg spiral}%
\label{voderberg}%
\end{figure}

One of the most beautiful spirals is the first ever published in the literature: the Voderberg spiral from
1936~\cite{voderberg} (see  Fig. \ref{voderberg}). We marked each tile with an arc as shown in the figure to generate a
connected curve without an endpoint. Although the principle of nonperiodicity is the same as above, the endless
curve of arcs can be regarded to be more satisfying. (The main result of Voderberg was not the construction of a
spiral which was sort of a by-product. The stunning feature of his prototile was the ability to be enclosed by two
of its copies.)
\begin{figure}[!htb]
\centering
\includegraphics[width=0.55\columnwidth]{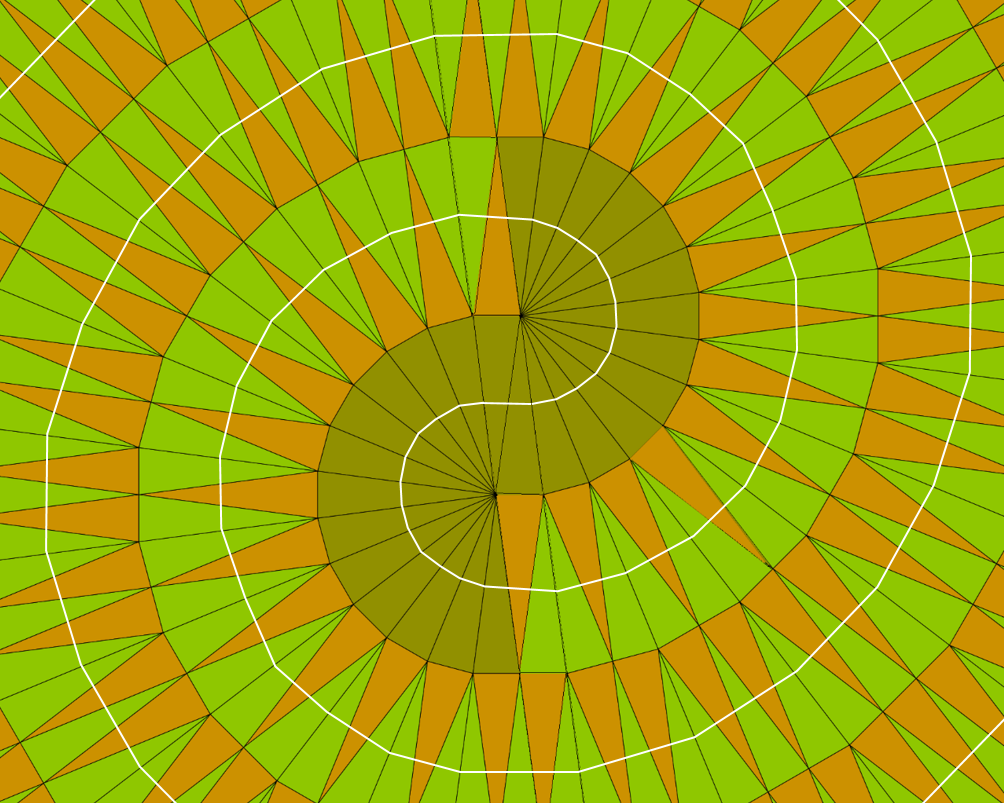}%
\caption{The triangle spiral}%
\label{triangle}%
\end{figure}

Finally, the prize for the most simple aperiodic prototile (following rule R1) goes to an isosceles triangle with
angle at the tip dividing $\pi$ (here we chose $\pi/12$ as in the example tiling from~\cite{TP}) and the marked arc
halving the legs (see Fig. \ref{triangle}). This shockingly simple tile together with rule R1 fulfills the conditions
of Lemma \ref{lem:curve} and therefore (as the other examples) only generates a nonperiodic (i.e.~spiral) tiling.

\section{Other types of rules}
Many different sorts of rules (not only matching rules) were discussed in the literature, so one can group them in this way:

\begin{description}
\item[Rule type 1] Matching rules for adjacent tiles (no colors required)~\cite{penta}

\item[Rule type 2] Colors required for matching adjacent tiles ~\cite{hpt}

\item[Rule type 3] Non-adjacent but pairwise matching within a neighborhood~\cite{Aaht}

\item[Rule type 4] Configurations of permitted adjacent tiles from an atlas~\cite{open}

\item[Rule type 5] Checks with adjacent tiles following a sequence or tree  (\cite{dendrites} or above cases)
\item[Rule type 6]  A given number of translation classes must be found in the tiling (see below).
\item[Rule type 7]  A certain patch of tiles must exist in the tiling (\cite{grow} and see below).
\end{description}

We can give a further distinction: For rule types 1--4 it is possible to perform a rule check for each finite region of the tiling locally (i.e.,~just by inspecting the set of connected parts each containing a small number of tiles). Usually these checks are made with adjacent tiles but there are important examples where non-adjacent pairs of tiles have to be inspected (rule type 3).  

For rule types 5--7 a similar check is also possible, however must be performed during a defined process of
generating the tiling: In~\cite{dendrites} a tree-like structure is built, in our approach in the above sections a
continuous curve or sequence of tiles should be followed. This group of rule types is working with a starting set or
seed, which can be an arbitrarily chosen first tile or (as we will see) a certain combination of tiles that is
demanded to occur in the tiling by the given rule. 

Rule types 6 and 7 check for the existence of certain tiles or angular positions of tiles. Also in these cases a seed set would be the simplest form of check to make sure that the demanded objects can be found in the tiling without a lengthy searching process. It should be pointed out that the use of a seed tile (or set of tiles) is a quite demanding restriction. Note that a seed set can be regarded as a ``rule of existence'' which principally cannot be checked locally if the full tiling is given. At first glance one might think that a single and arbitrary seed tile (as required in rule type 5) should not be a strong restriction, but it cannot be neglected that such a seed tile is the nucleus for the spiral structure in all our examples so far. 

\begin{figure}[!htb]
\centering
\includegraphics[width=0.6\columnwidth]{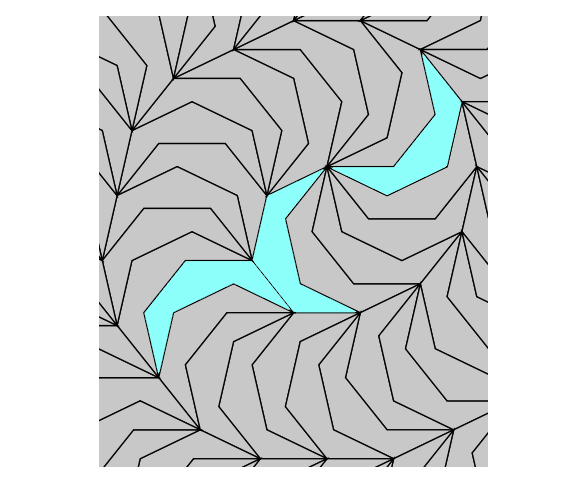}%
\caption{These three highlighted tiles force a spiral}%
\label{Three}%
\end{figure}

\begin{figure}[!htb]
\centering
\includegraphics[width=0.5\columnwidth]{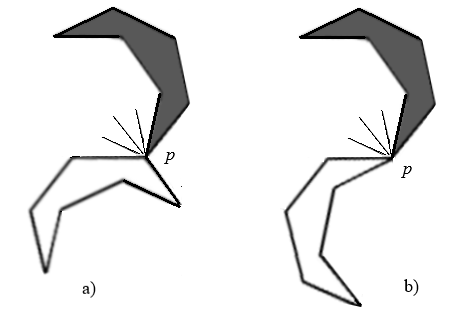}%
\caption{These gaps can only be filled in a unique way}%
\label{Filling}%
\end{figure}

Also in the following case with rule type 7 the seed is the key element for the final structure, here working without any decoration or modification of the prototile:

\begin{prop} \label{prop:2}
A patch of three tiles formed as the highlighted tiles in Figure \ref{Three} forces a nonperiodic tiling with prototile $T$ without any decoration or further matching rule.
\end{prop}
For the proof of this proposition we need the following lemmas.
\begin{lem} \label{lem:2}
Given two copies of $T$ positioned as shown in Fig. \ref{Filling}, in a monohedral tiling the gap of $4\pi/7$ can only be filled by copies of the black tile rotated around $p$, the vertex shared by both tiles. 
\end{lem} 

\begin{proof}
At point $p$ the only way to fill the open angle of $4\pi/7$ is by placing four additional tiles all sharing $p$. For the new tile closest to the white tile, the convex side cannot be placed to the convex side of the white tile, since then the new tile must intersect the black one. Therefore, there is only one way to place it: The concave side of the new tile must touch the convex side of the white tile. With the same argument the other three tiles are uniquely placed to fill the gap.
\end{proof}

The proof shows that the same holds for analogue constellations with gaps of angle $3\pi/7$, $2\pi/7$ or $\pi/7$. 

\begin{figure}[!htb]
\centering
\includegraphics[width=1.0\columnwidth]{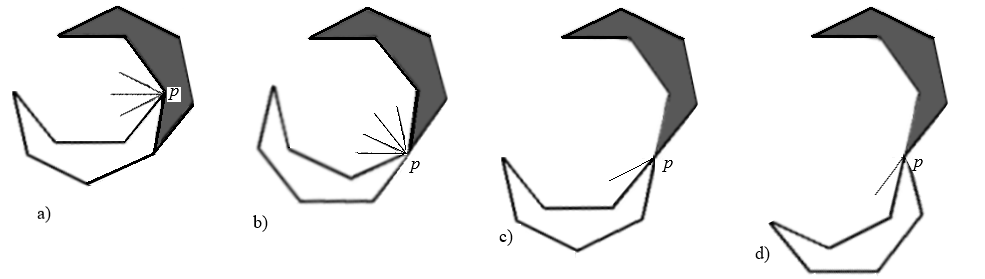}%
\caption{These constellations are impossible}%
\label{ForbiddenABC}%
\end{figure}

\begin{lem} \label{lem_forb}
The constellations shown in Figure \ref{ForbiddenABC} are impossible (in the sense that these configurations cannot occur in a monohedral tiling). 
\end{lem} 

\begin{proof} 
Case a: At point $p$ the only way to fill the open angle of $4\pi/7$ is by placing four additional tiles all sharing $p$. However, this is impossible, since one of these added tiles sharing an edge with the dark tile cannot be placed without overlapping with the white tile. Therefore, this constellation and its reflected copies cannot exist in a monohedral tiling. 

 Case b: To fill the gap at point $p$ with five tiles will not work, since these tiles must intersect each other. The only way to fill the gap with a single tile and avoiding constellation a) is shown in Fig. \ref{Forb2}.
Here a new gap with angle $4\pi/7$ occurs which cannot be filled without intersecting one of the tiles.

\begin{figure}[!htb]
\centering
\includegraphics[width=0.3\columnwidth]{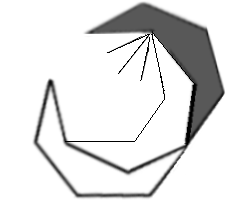}%
\caption{An impossible way to fill the gap in case b) }%
\label{Forb2}%
\end{figure}
\;\\
Case c: By filling the gap at point $p$, we must create a new constellation which is equivalent to b).
\;\\
 Case d: (Same argument) by filling the gap, we generate a constellation that was forbidden in c).
\end{proof}

\begin{proof}[Proof of Proposition \ref{prop:2}]
\begin{figure}[!htb]
\centering
\includegraphics[width=0.5\columnwidth]{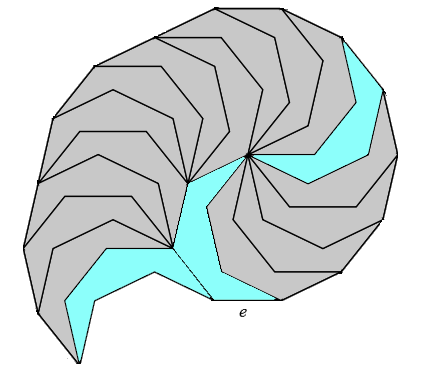}%
\caption{After filling the gaps between tiles }%
\label{patch}%
\end{figure}
With the same arguments as in Lemma \ref{lem:2} we can fill the gaps between the three tiles in a unique way and thus end up with the patch shown in Fig. \ref{patch}. The next question is: What happens at edge $e$? Is there a unique continuation of the spiral? To answer this question, we may have a look at Fig. \ref{Edge_e3}. The first and second cases generate an open gap at point $p$ at which the only way to fill this angle leads to a forbidden constellation shown in Lemma \ref{lem_forb}, case $a$.  
\begin{figure}[!htb]
\centering
\includegraphics[width=1.03 \columnwidth]{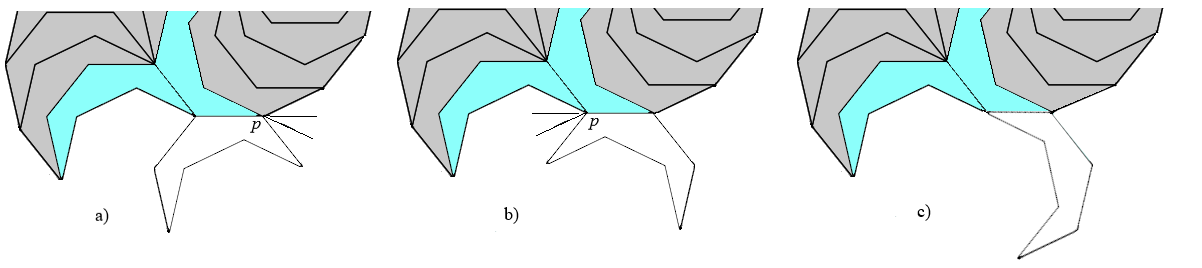}%
\caption{Forbidden constellations }%
\label{Edge_e3}%
\end{figure}
The rightmost case in Fig. \ref{Edge_e3} is equivalent to case $d$ from the said lemma and thus also forbidden. Therefore, there is only one unique continuation here as shown in the spiral tiling from Fig. \ref{spiral}. And the same arguments are valid for all following edges along the boundary of the growing patch. The inner obtuse angles at the patch's boundary are always $6\pi/7$ (e.g.,~at the points $p$), thus the situation at edge $e$ is principally the same for all following edges. So, we have shown that there is no other way to generate a monohedral tiling from the seed constellation of the three highlighted tiles.
\end{proof}

If we regard the patch $P$ from Fig. \ref{patch} as one tile, we come to another observation:

\begin{rem} \label{rem:4}
 Tile $T$ and patch $P$ (here considered as one tile) form an aperiodic tile set of two tiles without any matching rule. The only rule is here that both tiles should exist in the tiling (a simple variant of rule type 7).
\end{rem} 
\begin{proof}
If any tiling with translational symmetry is composed of the two tiles and at least one tile $P$ occurs in the tiling, then from Proposition \ref{prop:2} we learned that the existence of a patch $P$ must lead to the spiral tiling shown in Fig. \ref{spiral}. This contradicts the assumption that the tiling has translational symmetry.
\end{proof}

\begin{rem} \label{rem:5}
Other nonperiodic monohedral tilings are possible with tile $T$ and rule R1. It is not so important for this paper to show them, so it might be fun for the reader to explore them on his/her own. (They all have a spiral structure.)
\end{rem} 

At the beginning of this section we suggested another rule type (see type 6 in the above list) and conjecture that in any periodic tiling with prototile $T$ not all 14 translation classes will be present as it is the case in the spiral tiling of Fig. \ref{spiral}. So, one could set up a rule in form of a number of different translation classes that should (at least) be found in a tiling to enforce nonperiodicity. This could also be a topic for further studies.

\section{A non-spiral case}
One could come to the conclusion that this seed approach to nonperiodic tilings always results in tilings with a spiral structure. However, the spiral structure is nowhere used in the proofs, so it would be sufficient to generate a continuous, single, self-avoiding curve that meets every tile. We are going to demonstrate this with a prototile that leads to a different type of space-filling curve. 

 \begin{figure}[!htb]
\centering
\includegraphics[width=0.4 \columnwidth]{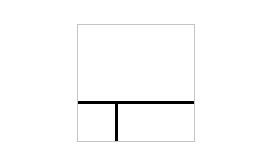}%
\caption{Another prototile with a decoration }%
\label{hilb1}%
\end{figure}
\begin{figure}[!htb]
\centering
\includegraphics[width=0.9 \columnwidth]{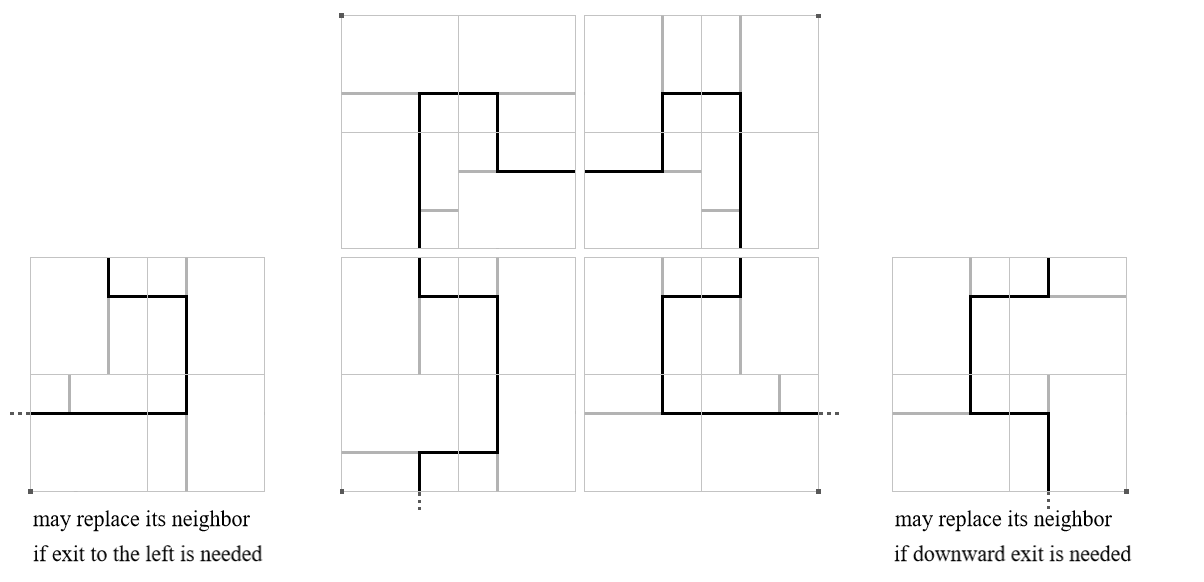}%
\caption{Elements for construction of a Hilbert polygon}%
\label{hilbert}%
\end{figure}

Regard the decorated prototile in Fig. \ref{hilb1}.
With this simple tile with three possible connections, we need a slightly different matching rule:
\noindent
\begin{itemize}[label={}]
\item \textbf{Matching rule R2}: The decoration of each additionally placed tile continuously connects to those of an already placed tile such that for each tile there remains an open connection.
\end{itemize} 
Obeying this rule, we again generate a sequence of tiles and the union of decorations must contain a single
unlimited and self-avoiding curve. It would be easy to construct a tiling with a spiral decoration. However, our
goal here is to end up  with a tiling, in which the decorations all are connected and show a structure that contains
a Hilbert polygon\footnote{The correct term would be ``polygonal chain''. Just for conciseness we use ``polygon''
here.}\cite{planefilling}. Consider the constructions in Fig. \ref{hilbert}.
Four copies of the prototile are combined to quadruples to form the characteristic parts of a Hilbert polygon. (The path for the resulting polygon is highlighted by black lines which is done just for the reader's convenience. This should not mean that several different decorations are used. The coloring indicates that dead ends can be ignored here as long as we have an ongoing curve through the tiling.)

 The four quadruples in the center of Fig. \ref{hilbert} form the typical section within this polygon, and the two quadruples to the left and to the right are candidates to replace the parts in the lower row of the quadruples (depending of the desired entry or exit of the polygon. 
  \begin{figure}[!htb]
\centering
\includegraphics[width=0.9 \columnwidth]{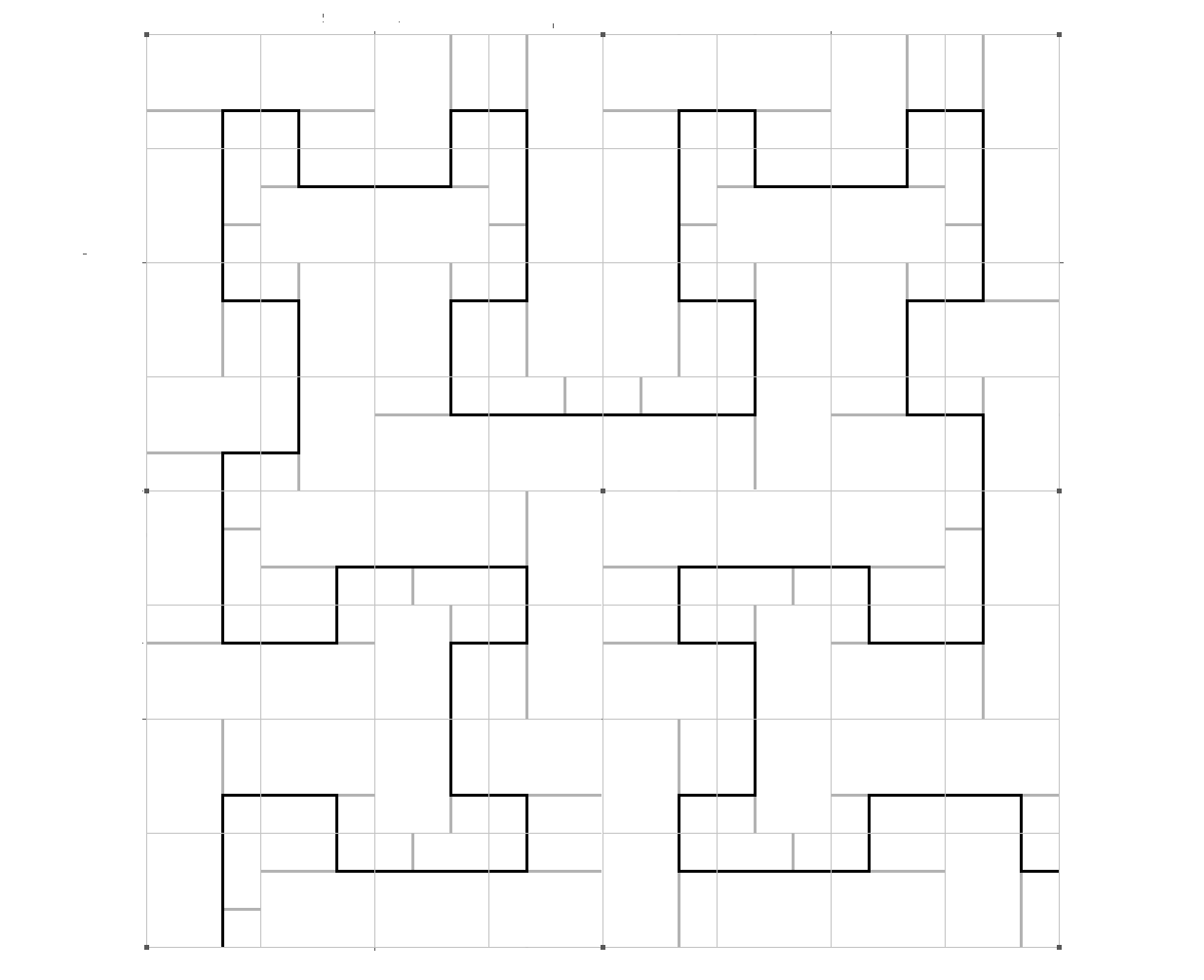}%
\caption{The union of decorations contains a Hilbert polygon }%
\label{hilbertBig}%
\end{figure}
Finally, Fig. \ref{hilbertBig} shows the bigger picture. At least to the author's knowledge, this is the first published construction of a Hilbert polygon by a decorated monohedral tiling.
Yet, two steps are missing: To demonstrate that this tiling fills the entire plane with a continuous polygonal curve and the proof that no self-intersection can occur on this path. Note that a Hilbert polygon is the plane-filling analogon to a Hilbert curve (continuously filling the unit square). \emph{Plane-filling} in this context means that for any point of the plane the distance to the polygon is not larger than a certain finite value. For the considered tiling, if it obeys rule R2, we know that the curve composed of the decorations on each tile will reach every tile. 
 \begin{figure}[!htb]
\centering
\includegraphics[width=1.05 \columnwidth]{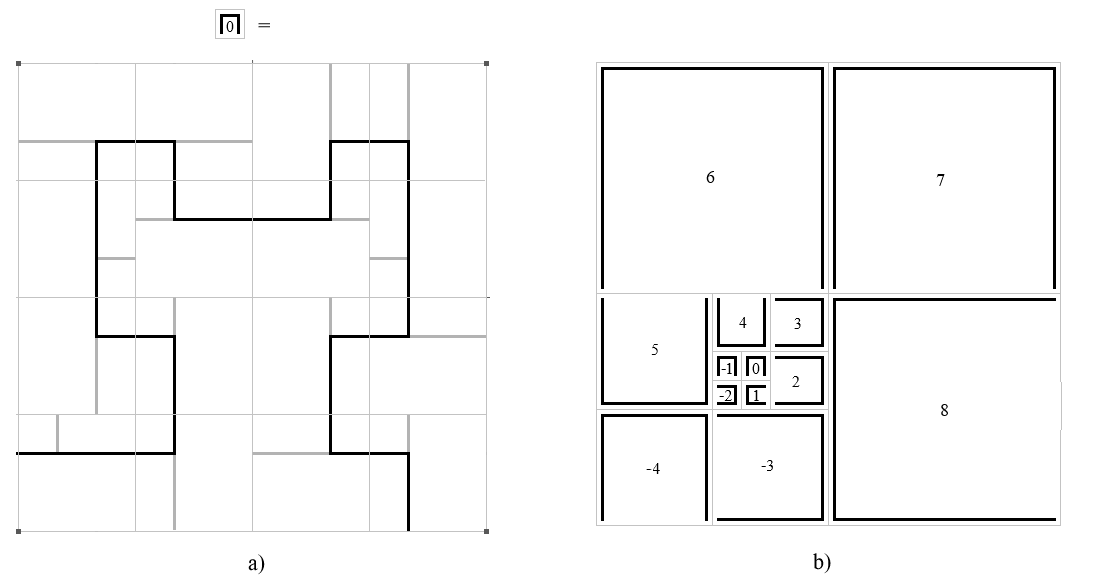}%
\caption{A schematic demonstration how the tiling is growing in a non-spiral manner in two directions}%
\label{hilbertGrow}%
\end{figure}

To make sure that R2 can be fulfilled, we demonstrate the growth scheme for the Hilbert tiling in Fig. \ref{hilbertGrow}.  On the left part of the Figure we describe a more abstract notation for a 16-tuple
of tiles, marked with number 0. Together with the other 16-tuples marked with 1, -1 and -2 (in the right half) they represent the four quarters of Fig. \ref{hilbertBig} and the starting cell of the plane tiling. The positive and negative signs in this numbering indicate that the curve is growing in two directions.

 The next bigger part (marked with 2) is a composition of 0, -1, -2, and 1, rotated by $90^{\circ}$ to the right. Only at the two corners of region 2 where the polygonal curve enters and exits, different variants of the quadruples can be inserted to force the curve into the appropriate direction (as explained in Fig. \ref{hilbert}). It is quite usual for a Hilbert curve that the connections to the neighbor regions have to be adapted due to the way the curve is running, however this is done in a well-defined manner, since there is only one appropriate choice in each case. Region 5 then is composed of the smaller regions -2 to 4 and again a rotation and an adaption of the connecting tiles has to be done. In this way the growth of the tiling ad infinitum can be described, which demonstrates that the full plane is covered.

 \begin{figure}[!htb]
\centering
\includegraphics[width=0.9 \columnwidth]{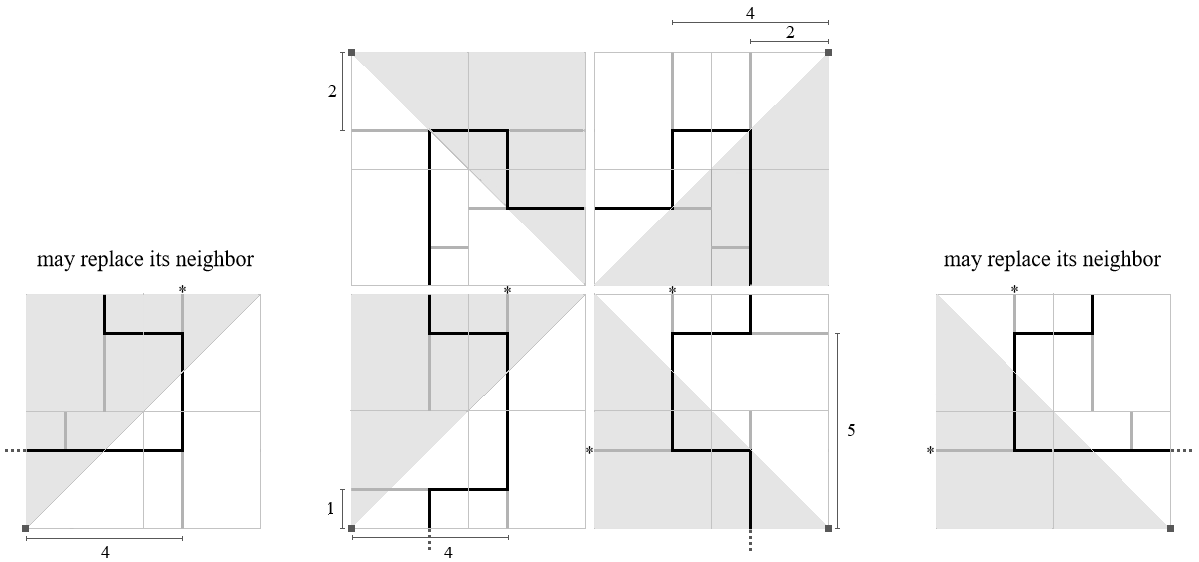}%
\caption{A coloring scheme to show that no self-intersections or loops via dead ends can occur }%
\label{hilbertColor}%
\end{figure}
Finally, we have to make sure that no extra loops or self-intersections can occur by the third connection that comes with each new tile. We can easily prove this by a coloring scheme shown in Fig. \ref{hilbertColor}.
 The corners of each 16-tuple (forming the basic cells of the tiling) are marked by tiny gray squares. The dead ends occurring at the boundaries can be characterized by their distance to these small squares. There are distances of length 1, 2, 4 and 5 (one square tile should have length 3). These distances w.r.t the closest small gray square are invariant under rotation of $90^{\circ}$ and under horizontal or vertical translation of 12 length units.

Observe the gray-white coloring of the triangles and note that this coloring also is invariant under rotation by multiples of $90^{\circ}$ and that triangles never share an edge with other triangles of the same color. Further note that the dead ends with distances 2, 4 and 5 only occur in white triangles, so they can never meet each other via a gray triangle. On the other hand, the dead ends with distance 1 only occur in gray triangles, so they never can meet other dead ends in white triangles. There is another type of dead ends (marked by asterisks) which occurs within the 16-tuple. They always are positioned at the same places and do not disturb each other. This shows that dead ends always remain dead ends in this tiling and thus the path of the generated Hilbert polygon is well-defined and self-avoiding.

Applying Proposition \ref{prop:1}, we can summarize that the presented tile from Fig. \ref{hilb1} together with rule R2 enforces a nonperiodic tiling. (This should not be misunderstood: The Hilbert-like tiling is not enforced by rule R2 but it is the most interesting tiling generated by this rule. Others, with simpler structure, are possible, but all of them must be nonperiodic.)

\section{Discussion}
In this study, the concept of plane-filling continuous curves (e.g.,~spirals) in form of the union of decorations on
a single prototile is linked to the well-known einstein problem. Although in its sharpest form this problem remains
unsolved, several aperiodic tile sets with size one and a single matching rule could be presented.  In addition, a new
type of rule was discussed which led to the same nonperiodic tilings. Also a plane-filling Hilbert polygon (a
discrete analogon to the Hilbert curve) is constructed by a monohedral decorated tiling. However, we should not
disguise that these approaches have one drawback: Despite the fact that the generated tilings have no translational
symmetry, the ``degree of nonperiodicity'' (to use an informal term) in the spiral tilings is not the same as in
Penrose tilings~\cite{penta} or other examples generated by known aperiodic tile sets. In fact, the periodic parts
of these tilings can become arbitrarily large in form of sectors, which is always the case in tilings with a spiral
structure. 

Nevertheless -- although the author is not an expert in quasicrystals or chemistry -- the concept of a starting set or
seed (or the demand for a certain constellation to exist in the investigated structure) seems to be a well-known
principle of crystal growth~\cite{grow} and applicable to such chemical structures where a crystal nucleus is needed
to create a new structure. Maybe a resulting sector structure as in our examples could also be useful if, for instance, certain optical effects are favorable. One should keep in mind that the first aperiodic tile sets were constructed
when nobody thought about practical applications in form of quasicrystals or in any other way. 

\section{Acknowlegdements}
I would like to thank Chaim Goodman-Strauss for valuable advices, and as well several hints from the journal reviewers were helpful to improve this paper.


\bibliographystyle{spmpsci}      
\bibliography{bibliog}   

%
%

\end{document}